\documentclass[a4paper,11pt]{amsart}

\usepackage{verbatim}
\usepackage{amsmath}
\usepackage{amssymb}
\usepackage{amsfonts}
\usepackage{amsthm}
\usepackage{mathrsfs}
\usepackage{tikz}
\usetikzlibrary{matrix,arrows}
\usepackage[all]{xy}
\usepackage{mathrsfs}
\usepackage{enumerate}
\usepackage{hyperref}
\usepackage[utf8]{inputenc}

\theoremstyle{plain}

\newtheorem{theorem}{Theorem}
\newtheorem{lemma}[theorem]{Lemma}
\newtheorem{proposition}[theorem]{Proposition}
\newtheorem{corollary}[theorem]{Corollary}
\newtheorem{question}[theorem]{Question}

\newtheorem*{theorem1}{Theorem}
\newtheorem*{corollary1}{Corollary}

\theoremstyle{definition}
\newtheorem{definition}[theorem]{Definition}

\theoremstyle{remark}

\newcommand{\lb}{\lbrace}
\newcommand{\rb}{\rbrace}

\title{Almost quasi-isometries and more non-exact groups}
\date{June 2015}
\author{Martin Finn-Sell}
\address{Martin Finn-Sell, Fakultät für Mathematik, Oskar-Morgenstern-Platz 1, 1090 Wien, Österreich.}
\email{martin.finn-sell@univie.ac.at}

\begin{document}

\begin{abstract}
We study permanence results for almost quasi-isometries, the maps arising from the Gromov construction of finitely generated random groups that contain expanders (and hence that are not $C^{*}$-exact). We show that the image of a sequence of finite graphs of large girth and controlled vertex degrees under an almost quasi-isometry does not have Guoliang Yu’s property A. We use this result to broaden the application of Gromov’s techniques to sequences of graphs which are not necessarily expanders. Thus, we obtain more examples of finitely generated non-$C^{*}$-exact groups.
\end{abstract}

\maketitle
\bibliographystyle{plain}

\section{Introduction}

For finitely generated discrete groups, $C^{\ast}$-exactness, or its geometric counterpart property A, provides a deep connection between operator algebraic and coarse geometric properties of a group. When a group is known to be $C^{\ast}$-exact many natural problems concerning that group, either from topology or algebra, have solutions: for example, in topology the Novikov conjecture concerning homotopy invariance of higher signatures is known to hold for $C^{\ast}$-exact groups (this requires techniques from both these fields, see for instance the introductory texts \cite{MR1728880,MR866507}). $C^{\ast}$-exactness also plays an important role for operator algebraists, as crossed product constructions provide concrete ways to produce new examples of $C^{\ast}$-algebras from old (see \cite{MR2391387} for a good introduction) and can help in classification of $C^{\ast}$-algebras through approximation properties and concrete understanding of tensor product norms as is outlined in \cite{MR1403994}. Studying pathological examples leads to new insights, and this places major importance on constructions of non-$C^{\ast}$-exact groups.

There are currently multiple methods of constructing $C^{*}$-non-exact groups that rely on the same basic idea, namely begin with the equivalence, for countable groups, between property A (a suitable notion of amenability for spaces) and $C^{\ast}$-exactness. We then find spaces that do not have property A and embed them into groups in a appropriate way. The resulting groups are not $C^{*}$-exact as property A passes to subspaces. 

More concretely, let $X$ be a graph without property A and $\mathcal{A}$ be a finite alphabet. We consider labellings $\mathcal{L}$ of the edges of $X$ by letters from $\mathcal{A} \cup \mathcal{A}^{-1}$. From such a labelling, we define a group $\Gamma$ via a presentation given by generators $\mathcal{A}$ subject to relators $\mathcal{R}$. The relators are determined by the words in the free group on the alphabet $\mathcal{A}$ that appear as a label of all possible cycles in $X$.

In order to guarantee that the group $\Gamma$ is not the trivial group one asks for \textit{small cancellation} conditions on the labelling $\mathcal{L}$. We also obtain, from such small cancellation conditions, the existence of a natural map $f : X \rightarrow \Gamma$ obtained by picking a basepoint in $X$ (mapped to the identity in $\Gamma$) and reading labelled paths. The small cancellation conditions give some metric control on this map - we then use the fact that $X$ does not have property A to deduce that $\Gamma$ also does not have property A - and thus is not $C^{\ast}$-exact.

The strength of the small cancellation condition one imposes then determine the properties of the map $f$ as well as the group theoretic properties of $\Gamma$. The main feature of small cancellation conditions is to ensure that relators in the presentation do not overlap too much: this is measured combinatorially using conditions on lengths of \textit{pieces} (common subwords of distinct relators). 

Possible small cancellation conditions one considers in this situation range from \textit{graphical} small cancellation (this is similar to the traditional small cancellation condition as in \cite{MR1812024} except that arbitrary graphs are labelled rather than just cycle graphs) to the more general \textit{geometric} small cancellation - each of these conditions relies on a different notion of piece. For references that contain more information on the specifics of these conditions, see \cite{exrangrps,DG-graphical}.

As we are interested in building non-trivial groups with particularly controlled embeddings of a fixed graph $X$ it is important to know that labellings of $X$ with small cancellation conditions exist. The techniques currently known that produce non-$C^{\ast}$-exact groups rely on probabilistic methods to show this existence. 

Indeed, probabilistic methods are good for providing statements concerning existence (that is, labellings will have the required small cancellation condition with probability greater than 0), but can also be used to show genericity (labellings will have the required small cancellation condition with overwhelming probability). In this context, we are looking for existence, but genericity is particularly interesting as it says something more about the ubiquity of a given property of a group in that specific random model.  

Recall that the \textit{girth} of a finite graph is the length of its shortest simple non-trivial cycle. A sequence of finite graphs $\lbrace X_{i} \rbrace_{i \in \mathbb{N}}$ has large girth if $girth(X_{i})$ tends to infinity in $i$. It was shown in \cite{MR2831267} that the spaces $X$ constructed from sequences of large girth do not have property A.

Currently, the graphs $X$ for which a small cancellation labelling does exist are obtained from sequences of finite graphs $\lbrace X_{i}\rbrace_{i\in \mathbb{N}}$ with large girth, vertex degree bounded between $3$ and $D$ and the ratio of diameter over girth uniformly bounded above. 

In this case, Osajda \cite{Osajda-nonexact}, using a Lovasz' local lemma variant, showed the existence of labellings that satisfy a graphical small cancellation condition by building on earlier work in \cite{MR3346927} and the ideas of \cite{MR1945373}. This condition yields embeddings $f: X \rightarrow \Gamma$ that are either isometric or coarse (depending how flexible one is in the construction), but only provides existence (as a consequence of technique in probability theory that was used). As the map $f$ is definitely a coarse embedding in this case, it is easy to see that these groups will not be $C^{\ast}$-exact using the forbidden subspace type argument outlined above.

On the other hand, Gromov in \cite{MR1978492} (see also Arzhantseva and Delzant \cite{exrangrps} for a comprehensive treatment) produces labellings that satisfy the geometric small cancellation condition, and this condition produces generic labellings (so certainly it provides existence) that come with many other desirable features.  Gromov's technique also provides a different type of embedding $f: X \rightarrow \Gamma$, that is central to the themes in this paper, called an almost quasi-isometry: 

\begin{definition}(\cite{MR1978492}, \cite{exrangrps}).
Let $\lbrace X_{i}\rbrace_{i\in \mathbb{N}}$ be a large girth sequence of finite graphs and let $\lbrace Y_{i} \rbrace_{i}$ be a sequence of metric spaces. A family of maps $\varphi_{i}:X_{i} \rightarrow Y_{i}$ is called an \textit{almost quasi-isometric embedding} if there exist positive constants $a,k, (b_{i})_{i}$ such that the following holds:
\begin{enumerate}
\item $a d(x,y) - b_{i} \leq d(\varphi_{i}(x),\varphi_{i}(y)) \leq k d(x,y)$;
\item $b_{i}=o(girth(X_{i}))$.
\end{enumerate}
\end{definition}

These embeddings are not necessarily coarse embeddings or fibred coarse embeddings \cite{MR3116568}.

The main results of this paper concern the relationship between almost quasi-isometric embeddings and the forbidden subspace argument presented for violating property A. We remark also that almost quasi-isometric embeddings are examples of so-called \textit{weak} embeddings introduced by Higson, Lafforgue and Skandalis in section 8 of \cite{MR1911663} in their work on counterexamples to the Baum-Connes conjecture with coefficients. 

The operator algebraic techniques of Higson, Lafforgue and Skandalis can be applied to show that a group that contains an almost quasi-isometrically embedded large girth \textit{expander} is not $C^{\ast}$-exact. The key point is that something stronger is true - they show that the expansion of the sequence can be detected using K-theory of a certain algebra constructed from the expander and the group.

The main result in this paper is that the weak expansion type property that large girth sequences already possess is good enough to produce obstructions for property A. More precisely, we prove the following permanence result:

\begin{theorem1}
Let $\lbrace X_{i} \rbrace_{i \in \mathbb{N}}$ be a large girth sequence with vertex degree bounded between $3$ and $D$, and let $\lbrace \varphi_{i}: X_{i} \rightarrow Y_{i}\rbrace_{i}$ be an almost quasi-isometric embedding from $\lbrace X_{i} \rbrace_{i}$ to a sequence of finite metric spaces $\lbrace Y_{i} \rbrace_{i}$. Then the coarse disjoint union $Y=\sqcup_{i} Y_{i}$ does not have property A.
\end{theorem1}

An immediate corollary of this result is the following:

\begin{corollary1}
Let $\lbrace X_{i} \rbrace_{i \in \mathbb{N}}$ be a large girth sequence with vertex degree bounded between $3$ and $D$, let $\Gamma$ be a finitely generated discrete group and let $\lbrace \varphi_{i}: X_{i} \rightarrow \Gamma\rbrace_{i}$ be an almost quasi-isometric embedding. Then $\Gamma$ is not $C^{\ast}$-exact.
\end{corollary1}

This corollary provides an extension to the work of Gromov and produces yet more examples of non-$C^{*}$-exact groups by applying the techniques of \cite{MR1978492,exrangrps} to families of graphs that are \textit{not} expanders. The benefit is that these examples still have a chance to coarsely embed into Hilbert space, unlike those that contain an expanding sequence. 

We give two methods of proof of our result as both techniques introduce new ideas that are of independent interest. The first is operator algebraic and is based on both the work of Higson, Lafforgue and Skandalis \cite{MR1911663} as well as more recent ideas arising from a characterisation of property A in terms of ghost operators given by Roe and Willett in \cite{MR3146831}. The second is coarse geometric and exploits the negation of a version of \textit{uniform local amenability}, defined originally in \cite{MR3054308}, for families of uniformly discrete metric spaces of bounded geometry.

\section*{Acknowledgements}
This author is partially supported by the ERC grant of Prof. Goulnara Arzhantseva “ANALYTIC” no. 259527. He would like to thank Goulnara Arzhantseva for suggesting the question and Rufus Willett for reading some early remarks about continuous functional calculus.

\section{Ghost operators, scaled Laplacians and weak expanders}\label{sect:ghost}
In this section we develop ideas concerning weak expansion for large girth familes and describe the characterisation of property A in terms of ghost operators given by Roe and Willett in \cite{MR3146831}. Their result relies on applying the continuous functional calculus to certain operators related to the discrete Laplacian and a weak expansion property that gives control over the spectrum of these Laplacian type operators. Below we recall precisely the definitions we require to explain the ideas of Roe and Willett \cite{MR3146831} that we use in the proof of the main result.

\begin{definition}
Let $X$ be a metric space. We say $X$ has \textit{bounded geometry} if for every scale $R>0$ there is a natural number $N_{R}$ such that every ball of radius $R$ in $X$ has cardinality at most $N_{R}$. We say $X$ is \textit{uniformly discrete} if there is a uniform lower bound on the distance between distinct points. 
\end{definition}

Given a uniformly discrete metric space $X$ of bounded geometry, there is a way to associate a $C^{\ast}$-algebra to $X$ that bridges operator algebraic properties with coarse geometric properties. Let $\ell^{2}(X)$ be the complex Hilbert space spanned by delta functions $\delta_{x}$ for each point $x\in X$. Any bounded linear operator $T\in \mathfrak{B}(\ell^{2}(X))$ can be uniquely represented as a matrix $(T_{x,y})$ indexed by $X \times X$ where the entries are defined by $T_{x,y}=\langle T\delta_{x},\delta_{y} \rangle$. 

For $T \in \mathfrak{B}(\ell^{2}(X))$ we can define the propagation of $T$  by the formula:
\begin{equation*}
Propagation(T):= \sup \lbrace d(x,y) \mid T_{x,y} \not = 0\rbrace. 
\end{equation*}
The collection of bounded linear operators with finite propagation is denoted $\mathbb{C}[X]$ and clearly forms a $\ast$-subalgebra of $\mathfrak{B}(\ell^{2}(X))$. The closure of $\mathbb{C}[X]$ in the operator norm of $\ell^{2}(X)$ is called the \textit{uniform Roe algebra of $X$} and is denoted by $C^{\ast}_{u}(X)$.

We are interested in particular operators in the uniform Roe algebra of certain coarsely disconnected spaces: 
 
\begin{definition}\label{def:coarselydisconnected}
Let $\lbrace X_{i} \rbrace_{i \in \mathbb{N}}$ be a sequence of finite metric spaces that are uniformly discrete with bounded geometry uniformly in the index $i$, such that $\vert X_{i} \vert \rightarrow \infty$ in $i$. Then we can form the \textit{coarse disjoint union} $X$ with underlying set $\sqcup X_{i}$, metric $d$ given by the metric on each component and $d(X_{i},X_{j})\rightarrow \infty$ as $i+j\rightarrow \infty$. Any such metric is proper and unique up to coarse equivalence.
\end{definition}

Coarse disjoint unions appear in many places in the literature surrounding the coarse Baum-Connes assembly conjecture (See \cite{higsonpreprint, MR1911663, explg1, MR2568691} for a few instances of this connection), or more generally for producing interesting spaces without property A \cite{MR2831267} that also coarsely embed into Hilbert space \cite{MR2899681}.

\begin{definition}
Let $X$ be a metric space with bounded geometry. The \textit{Laplacian on scale $R$} is a bounded linear operator acting on $\ell^{2}(X)$ that is defined pointwise by the formula:
	\begin{equation*}
		(\Delta_{R}f)(x)=\sum_{(x,y)\in E_{R}} f(x)\delta_{x} - f(y)\delta_{y}
	\end{equation*}
	where $E_{R}$ is the $R$-neighbourhood of the diagonal in $X\times X$ determined by the metric on $X$ and $f \in \ell^{2}(X)$. A computation shows that
	\begin{equation*}
	  \langle \Delta_{R}f,f\rangle = \sum_{(x,y)\in E_{R}} \vert f(x) - f(y) \vert^{2}.
	\end{equation*}
	For a graph $X$ the Laplacian on scale $1$ is the unnormalised combinatorial Laplacian of $X$, and one should think of $E_{R}$ as a generalised notion of edge for a arbitrary metric space $X$.
\end{definition}

On a coarse disjoint union of finite metric spaces $X=\sqcup X_{i}$, we get a decomposition of the $\Delta_{R}$ into components $\Delta_{R,i}$ (the Laplacian on scale $R$ for the metric space $X_{i}$). This corresponds to the fact that in a coarse disjoint union metric the set $E_{R}$ for $X$ breaks into a union of the sets $E_{R,i}$, the $R$-neighbourhood of the diagonal in $X_{i}$, and a single finite set $F_{R}$ that collects together all the pairs from beginning of the sequence that belong to different $X_{i}$. This means that in order to study $\Delta_{R}$ for such a space, it is enough to understand the $\Delta_{R,i}$ (because the corresponding operators $\Delta_{R}$ and $\oplus_{i} \Delta_{R,i}$ differ by a compact operator).

With the operators $\Delta_{R,i}$ in hand we can formulate the notion of weak expansion from \cite{MR3073258,MR3146831}.

\begin{definition}A sequence of finite metric spaces $\lbrace X_{i} \rbrace_{i \in \mathbb{N}}$ with bounded geometry form a \textit{weak expander} if there exists a scale $R>0$ and a constant $\epsilon >0$ such that for every $S>0$ there exists $i_{S}$ such that for all $i>i_{S}$ and every non-zero vector $\eta_{i} \in \ell^{2}(X_{i})$ with $diam(supp(\eta_{i}))<S$, the Laplacian $\Delta_{R,i}$ satisfies:
\begin{equation*}
\langle \Delta_{R,i}\eta_{i}, \eta_{i} \rangle > \epsilon \Vert \eta_{i} \Vert
\end{equation*}
\end{definition}

This is a considerable generalisation of expansion and it is known that large girth sequences are \textit{weak expanders} (because they are not uniformly locally amenable, c.f \cite{MR3054308} Section 5). We observe that they are in fact slightly stronger than weak expanders.

\begin{lemma}\label{lem:LG+}
Let $\lb X_{i}\rb$ be a sequence of finite graphs with large girth and vertex degrees between $3$ and $D$. Then there exists $\epsilon>0$ such that for every $\eta \in \ell^{1}(X_{i})$ with $diam(supp(\eta))\leq \frac{1}{2}girth(X_{i})$ we have:
\begin{equation*}
\sum_{(x,y)\in E(X_{i})} \vert \eta(x) - \eta(y) \vert > \epsilon \Vert \eta \Vert_{1}.
\end{equation*}
\end{lemma}
\begin{proof}
Suppose $\eta$ is a positive vector (without loss of generality). Now find a family of nested subsets $F_{j}$ satisfying $\cup F_{j} = supp(\eta)$ and positive real numbers $a_{j}$ such that $\eta = \sum_{j} \frac{a_{j}}{\vert F_{j} \vert}\chi_{F_{j}}$. Now we compute:
\begin{eqnarray*}
\sum_{x\in X_{i}} \sum_{\substack{y\in X_{i}\\ (x,y)\in E}} \vert \eta(x) - \eta(y) \vert & = & \sum_{j} \frac{a_{j}}{\vert F_{j} \vert} \sum_{x\in X_{i}} \sum_{\substack{y\in X_{i}\\ (x,y)\in E}} \vert \chi_{F_{j}}(x) - \chi_{F_{j}}(y) \vert \\
& = & \sum_{j} \frac{2a_{i}}{\vert F_{j} \vert} \sum_{\substack{(x,y)\in E \\ x \in F_{j}, y\in \partial F_{j}}} \vert \chi_{F_{j}}(x) \vert \\
& \geq & \sum_{j}\frac{2a_{j}}{\vert F_{j} \vert}\vert \partial F_{j} \vert \geq \frac{2}{D-1} \Vert \eta \Vert_{1}.
\end{eqnarray*}
\end{proof}

A standard argument allows one to find an explicit $c:=c(\epsilon)$ such that for any $\xi \in \ell^{2}(X_{i})$ with $diam(supp \xi) \leq \frac{1}{2}girth(X_{i})$ we have $\langle \Delta \xi, \xi \rangle > c \Vert \xi \Vert$. The proof of Lemma 3.2 in \cite{MR3146831} does this calculation explicitly computing a constant $c$ from $\epsilon$ and in this case, for $\epsilon = \frac{2}{D-1}$, the constant $c=\frac{2}{(D-1)^{2}D}$.

In fact, the above can actually be seen as a reformulation of the weak expander property in general:

\begin{proposition}
Let $\lb X_{i} \rb$ be a uniformly discrete family with uniform bounded geometry. Then $\lb X_{i} \rb$ is a weak expander if and only if there is a $c>0$ and a sequence $S_{i}$ of positive increasing real numbers such that for all vectors $\eta \in \ell^{2}(X_{i})$ of $diam(supp(\eta))$ less than $S_{i}$ we have: 
\begin{equation*}
\langle \Delta_{i} \eta, \eta \rangle > c \Vert \eta \Vert.
\end{equation*} 
\end{proposition}
\begin{proof} This follows from the  proof of Lemma \ref{lem:LG+}, but using 
\begin{equation*}
S_{i}:= \sup \lb diam(supp(\eta)) | \eta \in \langle \Delta_{i}\eta,\eta\rangle > c \Vert \eta \Vert \rbrace
\end{equation*}
instead of girth in the estimates.
\end{proof}

One of the coarse geometric applications of expansion and weak expansion is to construct \textit{ghost operators}.

\begin{definition}\label{Def:Ghost}
Let $X$ be a metric space of bounded geometry and let $T$ be a bounded linear operator of finite propagation on $\ell^{2}(X)$. $T$ is a \textit{ghost operator} if $\forall \epsilon >0$ there exists a bounded subset $B \subset X\times X$ such that $\vert T_{xy} \vert \leq \epsilon$ for all $(x,y) \in (X\times X) \setminus B$.
\end{definition}

Every compact operator is a ghost operator in the above definition, By a recent result of Roe and Willett \cite{MR3146831} property A is equivalent to the fact that every ghost operator is a compact operator.

\begin{theorem}(Roe-Willett \cite{MR3146831})\label{thm:RW2}
Let $X$ be a metric space of bounded geometry. Then $X$ has property A if and only if every ghost operator in $C^{\ast}_{u}(X)$ is compact.
\end{theorem}

For the application we have in mind we also require a generalisation of Proposition 3.1 from \cite{MR3146831}, which proves that weak expanders admit ghost operators via an application of continuous functional calculus. We state this generalisation below the following definition.

\begin{definition}
Let $\lb X_{i} \rb$ be a sequence of finite graphs. Then a family of bounded linear operators $T_{i} \in \mathfrak{B}(\ell^{2}X_{i})$ are said to be asymptotically $\epsilon$-expanding if  for every $S>0$ there is an index $i_{S}$ such that for all $i\geq i_{S}$ and $\eta \in \ell^{2}(X_{i})$ with $diam(supp(\eta_{i}))<S$ we have $\langle T_{i}\eta, \eta \rangle > \epsilon \Vert \eta \Vert$.
\end{definition}

We remark that if a sequence of finite graphs is a weak expander then there are $R,\epsilon>0$ such that $\lbrace \Delta_{R,i} \rbrace $ is an asymptotically $\epsilon$-expanding family of propagation $R$ operators that are positive and self adjoint. This is the key observation that goes into the original version of the Proposition below.

\begin{proposition}(Roe-Willett \cite{MR3146831})\label{thm:RW+}
Let $\lbrace X_{i} \rbrace_{i \in \mathbb{N}}$ be a sequence of finite graphs with cardinality tending to infinity and vertex degrees uniformly bounded above by $D$ and below by $3$. Fix $R>0$, and let $\lbrace T_{i} \rbrace$ be an asymptotically $\epsilon$-expanding family of bounded linear operators that are positive, self adjoint and have uniform propagation at most $R$ such that the kernel of $T=\oplus_{i} T_{i}$ is infinite dimensional. Then by applying functional calculus to $T=\oplus_{i} T_{i}$ there is a non-compact ghost $f(T)$ in the uniform Roe algebra of $X=\sqcup X_{i}$.\qed
\end{proposition}

The proof of Proposition \ref{thm:RW+} is exactly the proof of Proposition 3.1 from \cite{MR3146831}.

We can now state and prove the first main result:

\begin{theorem}\label{thm:MR1}
Let $\lb X_{i} \rb_{i}$ be a sequence of finite graphs of large girth and vertex degrees bounded between $3$ and $D$, $X=\sqcup_{i} X_{i}$ be the coarse disjoint union of the $X_{i}$. Consider a family of finite metric spaces $\lb Y_{i} \rb_{i}$ with an almost quasi-isometric embedding $\lb \varphi_{i}: X_{i} \rightarrow Y_{i} \rb_{i}$. Then $Y = \sqcup Y_{i}$ does not have property A.
\end{theorem}
\begin{proof}
Our aim is to construct a family of operators that satisfy Proposition \ref{thm:RW+}. This would then produce a non-compact ghost operator, which in turn obstructs property A by appealing to Theorem \ref{thm:RW2}. We begin by observing that it is sufficient to check the case that the maps $\varphi_{i}$ are surjective, as we can restrict to the images to obtain surjections.

Let $\theta_{i}: \ell^{2}Y_{i} \rightarrow \ell^{2}X_{i}$ be the partial isometry defined by: 
\begin{equation*}
(\theta_{i}f)(x) = \frac{1}{\sqrt{\vert \varphi_{i}^{-1}(\varphi_{i}(x))\vert}} f(\varphi_{i}(x))
\end{equation*} 
for $f \in \ell^{2}Y$. Denote by $\theta := \oplus \theta_{i} : \ell^{2}Y \rightarrow \ell^{2}X$ the diagonal composition of these partial isometries. 

Let $\Delta_{i}$ be the graph Laplacian of $X_{i}$ (i.e $\Delta_{i}=\Delta_{1,i}$). Using $\theta_{i}$ we define operators $D_{i}$ using the formula $\theta_{i}^{*}\Delta_{i}\theta_{i}$. We will show that the family $\lbrace D_{i} \rbrace$ satisfies the hypotheses of Proposition \ref{thm:RW+}. An elementary calculation tells us that each $D_{i}$ has matrix entries that are $0$ outside of a $2k$-neighbourhood of the diagonal (this follows from the fact that $\Delta$ has propagation at most $2$) and so the family has finite propagation uniformly in $i$. It remains to control the spectrum of the $D_{i}$ appropriately. Let $D(\Delta)$ denote the orthogonal sum of the $D_{i}$ alternatively defined by the formula $D(\Delta)=\theta^{*}\Delta\theta$.

From the formula for each $D_{i}$ it is enough understand the size of the supports of the elements $\theta_{i}f$, where $f \in \ell^{2}Y_{i}$, in terms of the support of $f$: as if the support $\theta_{i}f$ can be chosen small enough then we can use Lemma \ref{lem:LG+} to conclude the result directly from the equation:
\begin{equation*}
\langle D(\Delta)f,f \rangle = \langle \theta^{*}\Delta\theta f,f \rangle = \langle \Delta\theta f,\theta f \rangle.
\end{equation*}
However, as $\lb \varphi_{i} \rb$ is an almost quasi-isometric embedding, we can conclude that for every $S>0$ and $i$ large enough the pre-image of a subset of $Y_{i}$ of diameter at most $S$ is a set of diameter at most $\frac{S+b_{i}}{a}$. Observing that a linear function of $b_{i}$ is still $o(girth(X_{i}))$  we can find $i_{S}$ such that $\frac{S+b_{i}}{a} \leq \frac{girth(X_{i})}{2}$. Applying Lemma \ref{lem:LG+} we obtain the required control of the spectrum of $D(\Delta)$, and it follows that there is a non-compact ghost operator in the uniform Roe algebra of $Y$.
\end{proof}

\section{Uniform local amenability with respect to a measure and weak expansion}\label{sect:ULA}
In this section we will utilise a different characterisation of property A to attack the initial question from a different, purely geometric angle. 

\begin{definition}\label{def:UULA}(Uniform (ULA($\mu$))
A family of uniformly discrete metric spaces with bounded geometry $\mathcal{X}=\lbrace X_{\alpha} \rbrace_{\alpha}$ is uniformly (uniformly locally amenable with respect to measures)\footnote{This naming convention is both well established in the literature and very unfortunate.} if for every $R,\epsilon>0$ there is an $S>0$ independent of $\alpha$ such that for every $\alpha$ and for every probability measure $\mu_{\alpha}$ on $X_{\alpha}$ there is a finite set $F_{\alpha}$ of diameter less than $S$ such that:
\begin{equation*}
\mu_{\alpha}(\partial_{R}F_{\alpha}) < \epsilon \mu_{\alpha}(F_{\alpha})
\end{equation*}
where $\partial_{R}F_{\alpha}:=\lbrace x_{\alpha} | 0 < d(x_{\alpha},F_{\alpha}) \leq R \rbrace$ is the $R$-boundary of $F_{\alpha}$ in $X_{\alpha}$. 
\end{definition}

This is the uniformisation of the property introduced in \cite{MR3054308} and is equivalent to the family of spaces having property A uniformly (this is clear from the results in \cite{MR3054308}). For a sequence of finite graphs $\mathcal{X} = \lbrace X_{i} \rbrace_{i \in \mathbb{N}}$, it is known also that the family $\mathcal{X}$ has property A uniformly if and only if the space $X=\sqcup X_{i}$ has property A - an in depth discussion of this kind of permanence for families can be found in \cite{EG-permanence}. We will make use of the negation of this property.

\begin{definition}
A family of finite uniformly discrete metric space $\mathcal{X}=\lbrace X_{\alpha} \rbrace_{\alpha}$ of bounded geometry that satisfies the negation of Definition \ref{def:UULA} for a set of measures $\underline{\mu}:=\lbrace \mu_{\alpha} \rbrace_{\alpha}$ will be called a \textit{$\underline{\mu}$-weak expander}.
\end{definition}

Our aim now will be to construct, from a given large girth sequence $\lbrace X_{i} \rbrace_{i}$ and a surjective almost quasi-isometry $\lbrace \varphi_{i}:X_{i} \rightarrow Y_{i}\rbrace_{i}$, a family of measures $\underline{\mu}=\lbrace \mu_{i} \rbrace_{i}$ with which the family of metric spaces $\lbrace Y_{i} \rbrace_{i}$ becomes a $\underline{\mu}-$weak expander.

\begin{theorem}\label{thm:MR2}
Let $\lbrace X_{i} \rbrace_{i \in \mathbb{N}}$ be a sequence of finite graphs with large girth and vertex degree bounded between $3$ and $D$ and let $\lbrace\varphi_{i}: X_{i} \rightarrow Y_{i}\rbrace_{i}$ be a surjective almost quasi-isometry to a sequence of finite metric spaces $\lbrace Y_{i} \rbrace_{i}$. Then there exists a family of measures $\underline{\mu}=\lbrace \mu_{i} \rbrace$ on $\lbrace Y_{i} \rbrace_{i}$ such that the family $\lbrace Y_{i} \rbrace_{i}$ is a $\underline{\mu}$-weak expander.
\end{theorem}
\begin{proof}It is enough to determine a scale $R$, a threshold $\epsilon$ and a probability measure $\mu_{i}$ on each $Y_{i}$. We define, for a subset $F_{i}$ of $Y_{i}$:
\begin{equation*}
 \mu_{i}(F_{i}) = \frac{\vert \varphi_{i}^{-1}(F_{i}) \vert}{\vert X_{i} \vert}.
\end{equation*}

Now we wish to show that for some $R>0$ and for every $S>0$, there is an index $i_{S}$ such that for all $i>i_{S}$ and $F_{i}\subset Y_{i}$ with diameter less than $S$ that we have $\mu_{i}(\partial_{R}(F_{i})) > \epsilon \mu_{i}(F_{i})$. This translates, after checking the definition of $\mu_{i}$, to the inequality $\vert \varphi_{i}^{-1}(\partial_{R}(F_{i})) \vert > \epsilon \vert \varphi_{i}^{-1}(F_{i}) \vert$. 

Observe that the preimages $\varphi_{i}^{-1}(F_{i})$ have diameter at most $\frac{1}{a}(S+b_{i})$ and so we find $i_{S}$ by taking the first natural number such that this number is less than $\frac{1}{2}girth(X_{i})$ - we know such an $i_{S}$ exists because the family $\lbrace X_{i} \rbrace_{i}$ has large girth.

Considering everything said above, it is enough to understand the relationship between the boundary of the preimage $\partial (\varphi_{i}^{-1}(F_{i}))$ and the preimage of the boundary $\varphi_{i}^{-1}(\partial_{R}(F_{i}))$ for each $R$: in particular if we take $i>i_{S}$ and $F_{i}$ of diameter less than $S$ (so that the diameter of $\varphi_{i}^{-1}(F_{i})$ is less than $\frac{1}{2}girth(X_{i})$), it will be enough to prove that $\partial (\varphi_{i}^{-1}(F_{i})) \subset \varphi_{i}^{-1}(\partial_{R}(F_{i}))$ as
\begin{eqnarray*}
\frac{\mu_{i}(\partial_{R}(F_{i}))}{\mu_{i}(F_{i})} & = & \frac{\vert \varphi_{i}^{-1}(\partial_{R}(F_{i})) \vert}{\vert \varphi_{i}^{-1}(F_{i}) \vert} \\ & \geq & \frac{\vert \partial (\varphi_{i}^{-1}(F_{i})) \vert}{\vert \varphi_{i}^{-1}(F_{i}) \vert} > \frac{1}{D-1}.
\end{eqnarray*}
Where the final inequality is a consequence of Lemma 3.3 from \cite{MR2831267} and the fact that the sets we are considering have diameter less than half the girth (and so they are finite trees). 
 
However, this is clear for $R=k$ as each $\varphi_{i}$ is $k$-Lipschitz. So the result holds for $\mu_{i}$ as above, with $R=k$ and $\epsilon = \frac{1}{D-1}$.
\end{proof}

\section{Concluding remarks and questions}
We begin with the main Corollary of Theorems \ref{thm:MR1} and \ref{thm:MR2}:
\begin{corollary}\label{cor:MT}
Let $\lb X_{i} \rb_{i}$ be a large girth sequence with degrees bounded between $3$ and $D$,  and let $\lbrace \varphi_{i}:X_{i} \rightarrow \Gamma \rbrace_{i}$ be an almost quasi-isometric embedding, where $\Gamma$ is a countable discrete group with equipped with a proper left invariant metric. Then $\Gamma$ is not $C^{*}$-exact.
\end{corollary}
\begin{proof}
The observation is that, so long as the images of the $\varphi_{i}$ are disjoint, we can factorise $\varphi_{i}$ as a surjective almost quasi-isometric embedding onto its image and a isometric embedding of that image into $\Gamma$. We can however always find some group element to act isometrically, by left translation, so that the images are disjoint. The result then follows from either Theorem \ref{thm:MR1} or Theorem \ref{thm:MR2}.
\end{proof}

The proofs from Sections \ref{sect:ghost} and \ref{sect:ULA} illustrate that even under very weak conditions on the embeddings, the coarse geometry of spaces that contain interesting families of finite graphs can be quite badly behaved. For a weak expander we know that there are both geometric and functional defintions that are equivalent. In the situation it is possible to construct different notions of scaled Laplacian using measures $\mu_{i}$ and one naturally asks:

\begin{question}
How are the ghost operators that arise from $D(\Delta)$ and other Laplacians constructed from the $\mu_{i}$ related?
\end{question}

Finally, we remark that combining with the techniques of \cite{MR1978492,exrangrps}, we arrive at new groups that are not exact. The class of graphs to which \cite{MR1978492,exrangrps} can be applied is now larger as a consequence of not requiring the expansion property (at least in practice: expanders are ubiquitous from the probabilistic point of view). We formulate this as a proposition below:

\begin{proposition}
Let $\lbrace X_{i} \rbrace_{i \in \mathbb{N}}$ be a family of large girth graphs that have diameter to girth ratio bounded by $C$ and that are $b$-thin (see Definition 5.3 from \cite{exrangrps}). Then the finitely generated group $\Gamma$ obtained through randomly labelling the sequence $\lbrace X_{i} \rbrace_{i \in \mathbb{N}}$ according to the Gromov model is not $C^{*}$-exact. \qed
\end{proposition}

This leaves the following broad question:

\begin{question}
What are the positive permanence properties of almost quasi-isometries? Can we, for instance, transfer a coarse embedding, wall structure or asymptotic coarse embedding through such maps?
\end{question}

A positive answer to one or all of the above would possibly allow for the construction of new groups with coarse embeddings into Hilbert space that are not $C^{*}$-exact.

\bibliography{ref.bib}

\begin{thebibliography}{10}

\bibitem{MR1945373}
Noga Alon, Jaros{\l}aw Grytczuk, Mariusz Ha{\l}uszczak, and Oliver Riordan.
\newblock Nonrepetitive colorings of graphs.
\newblock {\em Random Structures Algorithms}, 21(3-4):336--346, 2002.
\newblock Random structures and algorithms (Poznan, 2001).

\bibitem{exrangrps}
Goulnara Arzhantseva and Thomas Delzant.
\newblock Examples of random groups.
\newblock {\em Available on the authors' website}, 2008.

\bibitem{MR2899681}
Goulnara Arzhantseva, Erik Guentner, and J{\'a}n {\v{S}}pakula.
\newblock Coarse non-amenability and coarse embeddings.
\newblock {\em Geom. Funct. Anal.}, 22(1):22--36, 2012.

\bibitem{MR3346927}
Goulnara Arzhantseva and Damian Osajda.
\newblock Infinitely presented small cancellation groups have the {H}aagerup
  property.
\newblock {\em J. Topol. Anal.}, 7(3):389--406, 2015.

\bibitem{MR3054308}
Jacek Brodzki, Graham~A. Niblo, J{\'a}n {\v{S}}pakula, Rufus Willett, and Nick
  Wright.
\newblock Uniform local amenability.
\newblock {\em J. Noncommut. Geom.}, 7(2):583--603, 2013.

\bibitem{MR2391387}
Nathanial~P. Brown and Narutaka Ozawa.
\newblock {\em {$C^*$}-algebras and finite-dimensional approximations},
  volume~88 of {\em Graduate Studies in Mathematics}.
\newblock American Mathematical Society, Providence, RI, 2008.

\bibitem{MR3116568}
Xiaoman Chen, Qin Wang, and Guoliang Yu.
\newblock The maximal coarse {B}aum-{C}onnes conjecture for spaces which admit
  a fibred coarse embedding into {H}ilbert space.
\newblock {\em Adv. Math.}, 249:88--130, 2013.

\bibitem{MR1978492}
Misha. Gromov.
\newblock Random walk in random groups.
\newblock {\em Geom. Funct. Anal.}, 13(1):73--146, 2003.

\bibitem{DG-graphical}
Dominik Gruber.
\newblock Groups with graphical {C}(6) and {C}(7) small cancellation
  presentations.
\newblock {\em Trans. Amer. Math. Soc.}, 2014.

\bibitem{EG-permanence}
Erik. Guentner.
\newblock Permanence in {C}oarse {G}eometry.
\newblock {\em Recent Progress in General Topology III}.

\bibitem{higsonpreprint}
Nigel Higson.
\newblock Counterexamples to the coarse {B}aum-{C}onnes conjecture.
\newblock {\em Preprint}, 1999.

\bibitem{MR1911663}
Nigel Higson, Vincent Lafforgue, and Georges Skandalis.
\newblock Counterexamples to the {B}aum-{C}onnes conjecture.
\newblock {\em Geom. Funct. Anal.}, 12(2):330--354, 2002.

\bibitem{MR1403994}
Eberhard Kirchberg.
\newblock Exact {${\rm C}^*$}-algebras, tensor products, and the classification
  of purely infinite algebras.
\newblock In {\em Proceedings of the {I}nternational {C}ongress of
  {M}athematicians, {V}ol.\ 1, 2 ({Z}\"urich, 1994)}, pages 943--954.
  Birkh\"auser, Basel, 1995.

\bibitem{MR1812024}
Roger~C. Lyndon and Paul~E. Schupp.
\newblock {\em Combinatorial group theory}.
\newblock Classics in Mathematics. Springer-Verlag, Berlin, 2001.
\newblock Reprint of the 1977 edition.

\bibitem{Osajda-nonexact}
Damian Osajda.
\newblock Small cancellation labellings of some infinite graphs and
  applications.
\newblock {\em Preprint (http://arxiv.org/abs/1406.5015)}, 2014.

\bibitem{MR2568691}
Herv{\'e} Oyono-Oyono and Guoliang Yu.
\newblock {$K$}-theory for the maximal {R}oe algebra of certain expanders.
\newblock {\em J. Funct. Anal.}, 257(10):3239--3292, 2009.

\bibitem{MR3146831}
John Roe and Rufus Willett.
\newblock Ghostbusting and property {A}.
\newblock {\em J. Funct. Anal.}, 266(3):1674--1684, 2014.

\bibitem{MR866507}
Jonathan. Rosenberg.
\newblock {$C^\ast$}-algebras, positive scalar curvature and the {N}ovikov
  conjecture. {II}.
\newblock 123:341--374, 1986.

\bibitem{MR3073258}
Hiroki Sako.
\newblock A generalization of expander graphs and local reflexivity of uniform
  {R}oe algebras.
\newblock {\em J. Funct. Anal.}, 265(7):1367--1391, 2013.

\bibitem{MR2831267}
Rufus Willett.
\newblock Property {A} and graphs with large girth.
\newblock {\em J. Topol. Anal.}, 3(3):377--384, 2011.

\bibitem{explg1}
Rufus Willett and Guoliang Yu.
\newblock Higher index theory for certain expanders and {G}romov monster
  groups, {I}.
\newblock {\em Adv. Math.}, 229(3):1380--1416, 2012.

\bibitem{MR1728880}
Guoliang Yu.
\newblock The coarse {B}aum-{C}onnes conjecture for spaces which admit a
  uniform embedding into {H}ilbert space.
\newblock {\em Invent. Math.}, 139(1):201--240, 2000.

\end{thebibliography}

\end{document}